\documentclass{amsart}

\usepackage {amsfonts,amsmath,amsthm,amscd,amssymb}
\usepackage{srcltx}
\usepackage{hyperref}

\newcommand{\Z}{\mathbb{Z}}
\newcommand{\F}{\mathbb{F}}

\DeclareMathOperator{\Integer}{Int}
\providecommand{\Int}[1]{\Integer \left( #1 \right)}

\theoremstyle{plain}
\newtheorem{theorem}{Theorem}[section]
\newtheorem{conjecture}[theorem]{Conjecture}

\theoremstyle{definition}
\newtheorem{definition}[theorem]{Definition}
\newtheorem{example}[theorem]{Example}

\theoremstyle{remark}
\newtheorem{remark}[theorem]{Remark}

\numberwithin{equation}{theorem}
\numberwithin{table}{section}

\begin{document}

\title{Some conjectures and results about multizeta values for $\F_q[t]$}
\author{Jos\'e Alejandro {Lara Rodr\'iguez}}
\address{Universidad Aut\'onoma de Yucat\'an, M\'erida, Yucat\'an, M\'exico}
\email{lrodri@uady.mx, alara@math.arizona.edu }
\date{August 29, 2009}

\thanks{The author was supported by the Conacyt and the Universidad Aut\'onoma
de Yucat\'an}

\begin{abstract}
In this paper,  we explain several conjectures about how a product of two 
Carlitz-Goss zeta values can be expressed as a $\F_p$-linear combination of 
Thakur's multizeta values \cite{Thakur}, generalizing the $q=2$ case dealt
in \cite{Thakur_Multizeta08}. In contrast to the classical sum shuffle,
$\zeta(a)\zeta(b)=\zeta(a+b)+\zeta(a,b)+\zeta(b,a)$, the identities we 
get are much more involved. Hundreds of instances of these conjectures have
been proved and  we describe the proof method and the evidence. 

\end{abstract}

\maketitle

\setcounter{section}{-1}
\section{Introduction}

There are many recent works on multizeta values. We refer
to \cite{Waldschmidt} and references in there.

Thakur introduced multizeta related  to the arithmetic of function
fields, Drinfeld modules, $t$-motives and investigated the relations 
that they satisfied. We refer to \cite{A-T2}, \cite{Thakur} and
\cite{Thakur_Multizeta08} for details and we recall the definitions in
subsection ~\ref{notation}.

In this case, the identities are much more complicated than the classical
shuffle identities. In fact there are two types of identities, one with
$\F_p[t]$ coefficients and the other with $\F_p$ coefficients. In this paper, we
focus on the second type and,  in particular, explain several conjectures on 
how a product of two Carlitz zeta values can be expressed as $\F_p$-linear 
combination of Thakur's multizeta values \cite{Thakur}, generalizing the
$q=2$ case dealt in \cite{Thakur_Multizeta08}. In contrast to the classical
sum shuffle, $\zeta(a)\zeta(b)=\zeta(a+b)+\zeta(a,b)+\zeta(b,a)$, the 
identities we get are much involved. Hundreds of instances of these
conjectures have been proved and at the end we describe the proof method
and the evidence.

\subsection{Notation}
\label{notation}

\begin{align*}
\Z  & =  \{\mbox{integers} \}\\
\Z_+ & = \{\mbox{positive integers} \}\\
q & =  \mbox{a power of a prime } p\\
T,{t}& = \mbox{independent variables}\\
A & =  \F_q[t]\\
A_+ & =  \mbox{monics in } A\\
K & =  \F_q(t)\\
K_\infty & =  \F_q((1/t))=  \mbox{ completion of $K$ at $\infty$}\\
[n] & =  t^{q^n}-t\\
D_n & = \prod _{i= 0}^{n-1} (t^{q^n}-t^{q^i}) = [n][n-1]^q \dotsm [1]^{q^{n-1}
}\\
\ell_n & =  \prod _{i= 1} ^n (t-t^{q^i})= (-1)^n L_n = (-1)^n [n][n-1]\dotsm [1]
\\
\mbox{`even'} & =  \mbox{multiple of } q-1\\
\deg & =  \mbox{function assigning to } a\in A \mbox{ its degree in } t\\
\delta_{i,j} & = \mbox{the Kronecker delta} = \begin{cases}
                                                0 & \mbox{ for }i\ne j\\
                                                1 & \mbox{ for }i= j.\\
                                              \end{cases} \\
\Int{x}       & =   \begin{cases}
                    0 & \mbox{ if } x \mbox{ is not integer}\\
                    1 & \mbox{ if } x \mbox{ is integer}.\\
                    \end{cases} \\
\end{align*}

\subsection{Carlitz-Goss  zeta values}

For $s\in \Z$ and $d\in \Z_{\ge 0}$, write
\begin{align*}
S_d(s) := \sum _{ \substack{a\in A_+\\ \deg (a)=d}}\frac{1}{a^{s}} \in K
\end{align*}

For $s\in\Z$, the \emph{Carlitz-Goss zeta values} \cite{Goss96, Thakur} are defined as
\begin{align*}
\zeta_A(s):=\sum _{d=0}^\infty S_d(s) \in K_\infty.
\end{align*}

\subsection{Thakur multizeta values}

For $s_i\in \Z_+$, the \emph{Thakur multizeta values}  \cite{Thakur,
Thakur_Multizeta08} are defined by:

\begin{align*}
\zeta(s_1,\dotsc,s_r):= \sum _{d_1>\dotsb >d_r\ge 0} S_{d_1}(s_1) \dotsm
S_{d_r}(s_r)= \sum \frac{1}{a_1^{s_1}\dotsm a_r^{s_r} }\in K_\infty,
\end{align*}
where the second sum is over all $a_i\in A_+$ of degree $d_i$ such that 
$d_1>\dotsb >d_r\ge 0$. We say that this multizeta value has depth $r$ and
weight  $\sum s_i$.

\section{Relations between multizeta values} 

For $s_1,s_2\in \Z_+$ put
\begin{align*}
S_d(s_1,s_2)=\sum _{\substack {d=d_1>d_2\\a_i\in A_+}}\frac{1}{a_1^{s_1}a_2^{s_2}}
\end{align*}
where $d_i=\deg(a_i)$. For $a,b\in \Z_+$, we define
\begin{align*}
\Delta_d (a,b):=S_d(a)S_d(b)-S_d(a+b).
\end{align*}
The definition implies $\Delta_d(a,b)=\Delta_d(b,a)$.

We have:
\begin{conjecture}[Thakur, \cite{Thakur_Multizeta08}]\label{Delta(a,b)_as_linear_combination}
$\Delta_d(a,b)$ can be expressed as a linear combination  of $S_d$'s above:
\begin{align*}
\Delta_d(a,b) &=\sum c_j S_d(a_j,a+b-a_j),
\end{align*}
where $c_j\in \F_p^\times $ and $a_j$ are distinct positive integers. 
Furthermore, for each $a,b\in \Z_+$ the set $S(a,b)$ of pairs $(c_j,a_j)$ is
independent of $d$.
\end{conjecture}

Summing over the $d$'s we see  the product of two zetas as a linear  combination
of multizetas:
\begin{align*}
\zeta(a)\zeta(b)=\zeta(a+b)+\zeta(a,b)+\zeta(b,a)+\sum c_j\zeta(a_j,a+b-a_j).
\end{align*}

\textbf{Problem. } The problem is to give an algorithm to find $c_j,a_j$ given 
$(a,b)$ and $q$. In other words,  we want an algorithm to find $S(a,b)$ given
$(q,a,b)$. 

\begin{remark}

\begin{enumerate}
\item As we are mainly interested in conjecturing and proving how to get S(a,b)
given (q, a, b), by abuse of notation we write $\Delta(a,b)$ for
$\Delta_d(a,b)$, as this recipe is independent of d. In fact, recently Thakur
\cite{Thakur_shuffle09} has proved the conjecture as we state it, without
proving our conjectures in this paper.

\item A conjectural algorithm was given by Thakur \cite{Thakur_Multizeta08} for
$q=2$.

\item Since $S_d(p^n s) = S_d(s)^{p^n}$ and  $S_d(p^n s_1,p^n s_2) = S_d(s_1,s_2)^{p^n}$, 
  if $\Delta(a,b) =\sum c_j S_d(a_j,a+b-a_j)$, then $\Delta(a,b)^{p^n} = \Delta(p^n a,p^n b)$.
  Therefore, if $S(a,b)=\{(c_j,a_j)\}$,  then $S(p^n a,p^n b)= \{(c_j,p^n a_j)\}$.
  So we can restrict to $a$  and $b$ not divisible by $p$ without loss of generality.
\end{enumerate}
\end{remark}

We now give an illustrative example of several calculations, following
the method of \cite{Thakur_Multizeta08}, that we made using SAGE
\cite{Sage} and which form the basis of the conjectures below. 

\begin{example}\label{example_q5_a2}
Let  $q=5$, $a=2$ and $b=30$. By direct calculation,
\begin{multline*}
S_1(2)S_1(30)-S_1(32)= 3S_1(4,28)+3S_1(24,8) +4S_1(20,12)\\ +
S_1(12,20)+2S_1(8,24)+2S_1(28,4).
\end{multline*}
By conjecture \ref{Delta(a,b)_as_linear_combination}, we can express $\Delta(2,30)$ 
as a linear combination
\begin{multline*}
\Delta(2,30)= 3S_d(4,28)+3S_d(24,8) +4S_d(20,12)\\ +S_d(12,20)+2S_d(8,24)+2S_d(28,4).
\end{multline*}
Therefore, $S(2,30)=\{(3,4),(3,24),(4,20),(1,12),(2,8),(2,28)\}$. It may be
checked that this works for several $d$'s as is conjectured above.
In fact, it is proved below for all $d$.
\end{example}

\begin{definition}\label{main_definitions} Let $a\in \Z_+$.
\begin{enumerate}
\item  We set
    \begin{align*}
    r_a := (q-1)p^m,
    \end{align*}
    where $m$ is the smallest integer such that $a\le p^m$.
\item For $i,j$ put
    \begin{align*}
    \phi(i,j)&:=r_a-a-j(q-1)+ir_a,\\
    \phi(j)& :=\phi(0,j).
    \end{align*}
\item We define
    \begin{align*}
    j_{a, \max}:=\left \lfloor \frac{r_a-a}{q-1} \right \rfloor, 
    \end{align*}
    where $\lfloor x \rfloor$ is the largest integer not greater than $x$.

\item For $q$ prime, let $c_{a,j}\in \F_p ^\times$ be defined by:
    \begin{align*}
        c_{a,j}=\begin{cases}
        1 & \mbox{ if } j =0\\
        \left \lceil \frac{j(q-1)} {j_{a, \max} } \right \rceil ^{-1} 
        \binom{r_a-a}{j(q-1)} & 0 < j \le j_{a, \max}
    \end{cases}
    \end{align*}
    where $\lceil x \rceil $ is the smallest integer not less than $x$.

\item  For each $j,\;0\le j\le p-1$, let $\mu_j$ be the number of $j$'s in  the
$p$ expansion of $a-1$. Set
\begin{align*}
t_a = \prod _{j=0}^{p-2} (p-j)^{\mu_j}.
\end{align*}
\end{enumerate}
\end{definition}

The main conjecture is

\begin{conjecture}\label{main_conjecture} 

Let $a\in \Z_+$. Then 

\begin{enumerate}

\item  \label{r_a_for_any_q} The sets $S(a,b)$ can be found recursively by
\begin{align*}
S(a,b)= S(a,b-r_a)\cup T(a,b),
\end{align*}
where
\begin{itemize}
\item [a)] \label{t_a_for_any_q} $T(a,b)$ is a set of size $t_a$,
\item [b)] $T(a,b)$ has the form 
\begin{align*}
T(a,b)= \{(c_{j_{\ell} },b-\phi(j_{\ell}))\colon \ell = 0,\dotsc, t_a-1 \}
\end{align*}
for some $j_{\ell}$, $0\le j_{\ell}\le j_{a, \max}$.

\item [c)] The set $T_a$ of pairs $(c_{j_{\ell} },\phi(j_{\ell}))$ is independent of $b$.
\end{itemize}

\item \label{(1,phi0)_in_T_a}$(1,\phi(0)) = (1,r_a-a)\in T_a$.  Therefore if
$t_a=1$, then $T_a=\{(1,\phi(0))\}$.

\item \label{Ta_for_general_q} If there is no carry over base $p$ in the sum of
$j(q-1)$ and $\phi(j)$, then $(c,\phi(j)) \in T_a$ for some $c\in \F_p$.

\item For $r_a -q+2 \le b \le r_a$,
\begin{align*}
S(a,b) = \{(c_j,b-a_j)\colon (c_j,a_j)\in T_a\setminus \{(1,\phi(0))\} \}.
\end{align*}

\item \label{Ta_for_q_prime}If $q$ is prime, then $(c_{a,j},\phi(j)) \in T_a$ 
if and only if there is no carry over base $p$ in the sum of $j(q-1)$ and
$\phi(j)$.

\item \label{t_a=t_pma}If $T_a=\{(c_{j_{\ell}},\phi(j_{\ell}))\colon 0\le \ell
\le t_a-1\}$,  then $T_{p^ma}=\{(c_{j_{\ell} },p^m\phi(j_{\ell}))\colon 0\le
\ell\le t_a-1\}$.

\item Let us denote with $T^\theta(a,b)$  the set of $(c_j,a_j)\in T(a,b)$ such 
that $c_j=\theta$. Then the length of the sets $T^\theta(a,b)$ is independent of
$b$.
\end{enumerate}
\end{conjecture}

\begin{remark}
\begin{enumerate}
\item The conjecture about  $r_a$  improves the recursion estimate $(q-1)q^m$ 
given in \cite{Thakur_Multizeta08}.

\item Lucas theorem say that if $k =\sum k_i p^i$ and $m_j = \sum m_{ji}p^i$ are
 base $p$ expansions, then
\begin{align*}
\binom{k}{m_1,\dotsc,m_r}\equiv \prod \binom{k_i}{m_{1i}\dotsc, m_{ri} }\bmod p.
\end{align*}
Since we have that $\tbinom{a}{b}=0$ if $b>a$, we see as a corollary that the
multinomial coefficient displayed above is zero in $\F_p$ if there is a carry
over base $p$ in the sum $k=\sum m_i$.

In particular, if $q$ is prime we have that there is no carry over base $q$  in
the sum of $j(q-1)$ and $\phi(j)$ if
\begin{align*}
\binom{r_a-a}{j(q-1),\phi(j)}=\binom{r_a-a}{j(q-1)}\not\equiv 0 \bmod p.
\end{align*}
Therefore, Lucas theorem and  the conjecture about $T_a$ imply that when $q$ is
prime we have
    \begin{align*}
    T_a =\left\{ (c_{a,j},\phi(j))\colon \tbinom{r_a-a}{j(q-1)}\not\equiv 0\bmod q \right\},\\
    t_a = \# \{j\colon \tbinom{r_a-a}{j(q-1)}\not\equiv 0\bmod q,0\le j\le j_{a, \max} \}.
    \end{align*}

\item Since $c_{a,0}=1$ and $\tbinom{r_a-a}{0}=1$, part \ref{(1,phi0)_in_T_a} 
of \ref{main_conjecture} is consistent with part \ref{Ta_for_q_prime}.

\item Point \ref{t_a=t_pma} of \ref{main_conjecture} implies $t_a = t_{p^m a}$. 
This is consistent with our conjecture for the value of $t_a$. If $a-1=\sum a_i
p^i$ is the base $p$ expansion of $a-1$, then
\begin{align*}
p^m a -1 &= p^m -1 + \sum a_i p^{m+i}\\
& = (p-1)+(p-1)p+\dotsb + (p-1)p^{m-1}+ \sum a_i p^{m+i}\\
\end{align*}
is the base $p$ expansion of $p^m a -1$. This shows the equality of $t_a$ and
$t_{p^m a}$.

\item Give $q,a$ and $b$, a full conjecture for $S(a,b)$ includes
  \begin{enumerate}
  \item [a)] the recursion length $r_a$,
  \item [b)] the initial sets (or initial values) $S(a,b)$, $1\le b\le r_a$,
  \item [c)] the set $T_a$.
  \end{enumerate}
The main conjecture  does not give a distribution of signs for $T_a$ 
except when $q$ is prime, and does not take care of the initial values. In
section \ref{Conjectures_for_Delta(a,b)} we shall give a full description for
$q$ even and $a=2,3,4$ and also for $q$ odd and $a=2,3$. In section 
\ref{Full_conjecture_for_q=4} we shall give a full description  for $q=4$,
except for the initial values.

\end{enumerate}
\end{remark}

For $b>r_a$, let $b = r_a\sigma+b'$, $0< b' \le r_a$. Then
\begin{align*}
S(a,b)=S(a,b')\cup T(a,b-(\sigma-1)r_a)\cup \dotsb \cup T(a,b-r_a)\cup T(a,b),
\end{align*}
where each element of $T(a,b-ir_a)$ is of the form $(c_{j_{\ell}},b-\phi(j_{\ell})-ir_a)=(c_{j_{\ell}},\phi(i,j_{\ell}))$.
We see that  $\Delta(a,b)$ has the general form:
\begin{align}
\Delta(a,b) = \sum _{(\gamma,\alpha)\in S(a,b')} \gamma S_d(\alpha,\beta)+
\sum _{i=0}^{\sigma-1} \sum _{\substack{j\\ (c_j,\phi(j))\in T_a} } c_j S_d(b-\phi(i,j),a+\phi(i,j))
 \label{general_form_of_Delta(a,b)}
\end{align}
So $\Delta(a,b)$ has an  initial part (given by the initial values) and a regular (or recursive) part.

\begin{remark}
\begin{enumerate}

\item Thakur conjectured that all the $a+b-a_i$'s in conjecture \ref{Delta(a,b)_as_linear_combination} must be `even' \cite[5.3]{Thakur_Multizeta08}. Since $r_a$ is `even' we have that $a+\phi(i,j)=r_a-j(q-1)+ir_a$ is always even, which is consistent with this prediction.

\item If for every $a$ we know $T_a$ and we also know initial values $S(a,\beta)$ for $1\le \beta\le t_a$, we shall be able to compute $S(a,b)$ for any $b$. In \ref{main_conjecture} we define $T_a$ fully for $q$ prime.

\end{enumerate}
\end{remark}

\section{Full conjectures}
\label{Conjectures_for_Delta(a,b)}
In this section we present full conjectures, with initial values as well as the recursive recipe, for $\Delta(a,b)$ for small values of $a$ and also for some special large values of $a$ and $b$.

\subsubsection{\texorpdfstring{$\Delta(a,b)$}{Delta_a_b} for
\texorpdfstring{$a=2, 3, 4$}{a=2, 3, 4}, \texorpdfstring{$q$}{q} even}
When $q$ is a power of 2, we have $\F_2^\times=\{1\}$. To simplify notation
while describing $T_a$ instead of $(1,\phi(j))$ we shall use $\phi(j)$. From the
main conjecture \ref{main_conjecture}, it follows that if $q$ is prime and
$t_a=2$, then $T_a=\{(1, \phi(0)),(p-1,\phi(j_{a,\max}) \}$. But for example,
when $q=4$ and $a=7$, $T_7=\{\phi(0),\phi(3)\}$, but $3\ne 5 = j_{7,\max}$. We
conjecture that for any $q$ even we have $T_3 = \{\phi(0),\phi(j_{3,\max}) \}$.
Let $\phi,r_a$ and $j_{a,\max}$ as in definition \ref{main_definitions}. By
conjecture \ref{main_conjecture} and by this new conjecture we have

\begin{center}
\begin{tabular}{cccl} \hline
$a$ & $r_a$ & $t_a$ & $T_a$  \\ \hline
2 & $2(q-1)$ & 1 & $\phi(0)$  \\
3 & $4(q-1)$ & 2 &  $\phi(0),\phi(j_{3,\max})$\\
4 & $4(q-1)$ & 1 & $\phi(0)$  \\ \hline
\end{tabular}
\end{center}

\begin{conjecture}\label{Delta(q2n,2,b)}
Let $q$ be a power of $p=2$. Write $b=r_2\sigma+\beta$, $0<\beta\le r_2$. Then
\begin{align*}
\Delta(2,b)=\sum _{i=0}^{\sigma-1} S_d(b-\phi(i,0),2+\phi(i,0) )+
\Int{\tfrac{b}{q-1}}\tfrac{b}{q-1} S_d(2,b).
\end{align*}
\end{conjecture}

\begin{remark}
\begin{enumerate}

\item $\Int{x}$ is defined in section ~\ref{notation}.

\item Conjecture \ref{Delta(q2n,2,b)} is  a generalization of theorem 7 in \cite{Thakur_Multizeta08}.
\end{enumerate}
\end{remark}

\begin{conjecture}\label{Delta(q2,3,b)}
Let $q=2$.  Then
\begin{multline*}
\Delta(3,b)=\sum _{b-\phi(i,0)>3} S_d(b-\phi(i,0),3+\phi(i,0)) \\
+\sum _{b-\phi(i,j_{3, \max} )>3} S_d(b-\phi(i,j_{3, \max}),3+\phi(i,j_{3, \max}))\\
+ \sum _{i=1}^2 \Int{\tfrac{b-i}{r_3}} (S_d(2,b+1)+S_d(3,b)).
\end{multline*}
\end{conjecture}

\begin{conjecture}\label{Delta(q2n,3,b)}
Let $p=2$ and $q=p^n,n>1$. Then
\begin{multline*}
\Delta(3,b)=\sum _{b-\phi(i,0)>3} S_d(b-\phi(i,0), 3+\phi(i,0)) \\
+ \sum _{b-\phi(i,j_{3, \max})>3} S_d(b-\phi(i,j_{3, \max}),3+ \phi(i,j_{3, \max}))\\
+\Int{\tfrac{b+1}{q-1}}\tbinom{(b+1)/(q-1)+1}{2} S_d(2,b+1) \\
+ \Int{\tfrac{b}{q-1}}  (\tbinom{b/(q-1)+2}{2}-1) S_d(3,b).
\end{multline*}
\end{conjecture}

\begin{conjecture}\label{Delta(q2n,4,b)}
Let $q$ be a power of $p=2$. Then
\begin{multline*}
\Delta(4,b)=\sum _{b-\phi(i,0)>4} S_d(b-\phi(i,0),4+\phi(i,0)) \\
+ \Int{\tfrac{b-\max\{q-3,1\}}{r_4}} S_d(2,b+2) 
+ \Int{\tfrac{b-2q+3}{r_r}} S_d(3,b+1) \\
+ \sum _{i=1}^3 \Int{\frac{b-i(q-1)}{r_4}}S_d(4,b).
\end{multline*}
\end{conjecture}

\subsubsection{\texorpdfstring{$\Delta(a,b)$}{Delta_a_b} for
\texorpdfstring{$a=2, 3$}{a=2, 3}, \texorpdfstring{$q$}{q} odd} Next we give
full conjectures for $\Delta(a,b)$ for $a=2,3$ when $q$ is a power of an odd
prime. 
The problem becomes more complicated because me must take care of the different  coefficients, that is, we shall have to describe many  $S^\theta(a,b)$'s.

\begin{conjecture}\label{Delta(qpn,2,b)}
Let $p$ be an odd prime and let $q$ be a power of $p$. Then
\begin{multline*}
\Delta(2,b)=\sum _{j=0}^{p-1} (j+2)\sum _{b-\phi(i,p-1-j)>2} S_d(b-\phi(i,p-1-j),2 +\phi(i,p-1-j)) \\
+ \Int{\tfrac{b}{q-1}}\tfrac{b}{q-1} S_d(2,b).
\end{multline*}
\end{conjecture}

\begin{conjecture}\label{Delta(qpn,3,b)}
Let $q$ be a power of a prime $p>2$. Then
\begin{align*}
\begin{split}
\Delta(3,b) &= \sum _{b-\phi(i,0)>3} S_d(b-\phi(i,0),3+\phi(i,0)) \\
& +(1-\delta_{p,3}) \sum _{b-\phi(i,3) >3} S_d(b-\phi(i,3),3+\phi(i,3)) \\
& +\sum _{j=2} ^{\frac{p-3}{2}} \tbinom{j+1}{2}  \left (
\sum _{b-\phi(i,j+2)>3} S_d(b-\phi(i,j+2),3+\phi(i,j+2) ) \right .+ \\
& \left . \sum _{b-\phi(i,p+3-j-2)>3} S_d(b-\phi(i,p+3-j-2),3+\phi(i,p+3-j-2) )
\right) \\
& +(1-\delta_{p,3})\tbinom{\frac{p+1}{2}}{2} \sum _{b-\phi(i,\frac{p+3}{2})>3} S_d(b-\phi(i,\tfrac{p+3}{2}),3+\phi(i,\tfrac{p+3}{2}))\\ 
&  +\Int{\tfrac{b+1}{q-1}}\tbinom{(b+1)/(q-1)+1}{2} S_d(2,b+1) 
+ \Int{\tfrac{b}{q-1}}  (\tbinom{b/(q-1)+2}{2}-1) S_d(3,b).
\end{split}
\end{align*}
where $\delta_{i,j}$ denotes the Kronecker delta defined in section \ref{notation}.
\end{conjecture}

\subsubsection{Full conjectures for both indices large}

Here we have a conjecture for large indices, which are generalizations of
the formulas given in \cite[4.1.3]{Thakur_Multizeta08} 
\begin{conjecture}\label{Delta(a,b)_large_indices}
For general $q$ we have:
\begin{align}
\Delta(q^n,q^n-1) = 
- S_d(q^n,q^n-1).
\label{Delta(q^n,q^n-1)}
\end{align}
\begin{multline}
\Delta(q^n+1,q^n)  = \Int{\tfrac{2}{q}}S_d(2,2q^n-1) \\
- \sum _{j=1}^{\tfrac{q^n-1}{q-1}} S_d(3+(j-1)(q-1),2q^n-2-(j-1)(q-1) ).
\label{Delta(q^n+1,q^n)}
\end{multline}
\begin{multline}
\Delta(q^n-1,q^n+1)  = \\
- \sum _{j=1}^{\tfrac{q^n+q-2}{q-1}} S_d(2+(j-1)(q-1), 2q^n-2-(j-1)(q-1) ). 
\label{Delta(q^n-1,q^n+1)}
\end{multline}
\begin{multline}
\Delta(q^{n-1},q^n+1)  = \Int{\tfrac{2}{q}}S(2,q^{n-1}+q^n-1) \\
- \sum _{j=1}^{\tfrac{q^{n-1}-1}{q-1}} S_d(3+(j-1)(q-1), q^{n-1}+q^n-2+(j-1)(q-1) ). \label{Delta(q^{n-1},q^n+1)} 
\end{multline}

For $0\le i\le n$,
\begin{multline}
\Delta(q^n+1,q^n+1-q^i)  = \Int{\tfrac{2}{q}}S(2,2q^n-q^i) \\
- \sum _{j=1}^{\frac{q^n-q^i}{q-1}} S_d(3+(j-1)(q-1),2q^n-q^i-1-(j-1)(q-1)) \\
 + \sum _{j = \frac{q^n-q^i}{q-1}+1} ^{\frac{q^n-q^i}{q-1}}  S_d(3+(j-1)(q-1),2q^n-q^i-1-(j-1)(q-1) ).
\label{Delta(q^n+1,q^n+1-q^i)} 
\end{multline}
\end{conjecture}

\begin{remark}
\begin{enumerate}
\item Formula \eqref{Delta(q^n+1,q^n+1-q^i)} coincides with
\eqref{Delta(q^n+1,q^n)} when $i=0$. Furthermore,  when $p=2$, the sets
$S(q^n+1,q^n+1-q^i)$ are independent of $i$.
\item Conjecture \eqref{Delta(q^n+1,q^n+1-q^i)} generalizes the conjectured
phenomena \eqref{S(a,a-1)= S(a,a-4^j)} for $a=q^n+1$.
\end{enumerate}
\end{remark}

\begin{conjecture}
For $2\le m\le q$, we have
\begin{align*}
S_d(mq^i-1)
=\frac{\ell_{d+i} }{\ell_i \ell_d^{mq^i} }
=\frac{1}{\ell_d^{mq^i-1} }\frac{\ell_{d+i}}{\ell_i \ell_d}
\end{align*}
\end{conjecture}

We can deduce this for $m=2$ and $i=1$ from formula for $S_d(aq+b)$ in \cite[3.3.1]{Thakur_Multizeta08} and we have proved  for $m=i=2$ using the generating function \cite[3.2]{Thakur_Multizeta08} as well as verified it numerically for several low values of $m,i,q$ . This conjecture implies equation ~\eqref{Delta(q^n,q^n-1)} just as in Theorem 4 in \cite{Thakur_Multizeta08}.

\section{Full conjecture for \texorpdfstring{$q=4$}{q=4} except for initial
values}
\label{Full_conjecture_for_q=4}
In this section we present a prediction of how to write $\Delta(a,b)$ when $q=4$
as a sum of multizeta values if we know the initial values.

\begin{conjecture}\label{conjecture_for_q=4}
\begin{enumerate}
\item Given $a$ we compute $r_a$  and $t_a$ according to  \ref{main_definitions}.
\item Compute $T_a$ in the following way:
    \begin{enumerate}
    \item $\phi(j)\in T_a$ if there is no carry over base $p=2$ in the sum of $j(q-1)$ and $\phi(j)$ for $0\le j\le j_{a, \max}$.
    \item If there is carry over base $2$ in $j(q-1)$ and $\phi(j)$, let $\alpha_j$ be the number of 1's and 0's in the base 2 expansion of $j(q-1)$ and $\phi(j)$, if the total of 1's in the base 2 expansion of $j(q-1)$ and $\phi(j)$ equals 1 plus the number of 1's in the base 2 expansion of $r_a-a$. Let $M$ be the maximum of the $\alpha_j$'s. Then $\phi(j)\in T_a$ if $\alpha_j=M$.
    \end{enumerate}
\item Given $b>r_a$, write $b=r_a\sigma +b'$ where $0<b'\le r_a$. Then
\begin{align*}
\Delta(a,b)&= \sum _{a_i\in S(a,b')} S_d(a_i,a+b-a_i) +\sum _{i=0}^{\sigma-1} \sum _{ \substack{j\\ \phi(j)\in T_a}}S_d(b-\phi(i,j), a+\phi(i,j))
\end{align*}
\end{enumerate}
\end{conjecture}

Parallel to the $q=2$ case in \cite[4]{Thakur_Multizeta08}, we have the following conjecture:

\begin{conjecture}\label{S(a,a-1)= S(a,a-4^j)}
When $q=4$,
\begin{align*}
S(a,a-1)= S(a,a-4^j).
\end{align*}
\end{conjecture}

\section{Proofs}

As pointed out in \cite{Thakur_Multizeta08},  it is mechanical to 
verify multizeta identities of the type we are considering but one at a time, 
in other words, identities which do not have parameters.
Let us quickly recall the method.

The generating functions \cite[3.2]{Thakur_Multizeta08}  show that 
 $\ell_d^{k+1} S_d(k+1)$ and $\ell_d^k S_{<d}(k)$ are polynomials in $t^{q^d}$. We call $H_k$ and $G_k$ the polynomials in $K[T]$ such that 
\begin{align*}
H_k(t^{q^d}) &= \ell_d^{k} S_d(k), \\
G_k(t^{q^d}) &= \ell_d^{k} S_{<d}(k).
\end{align*}

To prove that
\begin{align*}\label{Delta=Delta'}
S_d(a)S_d(b)-S_d(a+b)= \sum _{(\gamma,\alpha)\in S(a,b)}\gamma S_d(\alpha,\beta),
\end{align*}
( where $a+b=\alpha+\beta$) for all $d$  
is equivalent to prove

\begin{align*}
H_a \left( T \right) \ H_b \left( T \right) - H_{a+b} \left( T \right)  & = 
\sum _{(\gamma,\alpha)\in S(a,b)}\gamma H_{\alpha} \left( T \right) \ G_{\beta} \left( T \right) \mbox{ in } \F_q(t)[T]
\end{align*}

By this method and using SAGE, we have proved the conjectures identities for 
$q=2, 4, 8, 3, 9, 5, 25, 13, 17$ and hundred values of $(a,b)$.

For details, more conjectures and identities we refer to the author's 
masters degree thesis \cite{Jalr} (also available by email). In
\cite{Thakur_shuffle09}, Thakur has described a more direct proof
method allowing to prove individual (without parameters) identities  faster.

{\bf Acknowledgments:} These results are 
 part of the Author's master's thesis 
 at the University of Arizona under the direction of Dinesh Thakur, to whom I am
very grateful for all his advice and unceasing encouragement over
the past year. His excellent guidance made these year of research a
pleasant experience. Thanks to Javier Diaz-Vargas for encouraging me to come
to
the University of Arizona. In addition, I want to express my gratitude
to the Universidad Aut\'onoma de Yucat\'an (Mexico) and the Consejo Nacional de
Ciencia y Tecnolog\'ia (Mexico) for their financial support.


\end{document}